\documentclass[12pt]{article}
\usepackage[cp866]{inputenc}
\usepackage{amssymb,amsmath,amsfonts}
\newtheorem{theorem}{{\bf Theorem}}[section]
\newtheorem{lemma}{{\bf Lemma}}[section]
\newtheorem{corollary}{{\bf Corollary}}[section]
\newtheorem{proposition}{{\bf Proposition}}[section]
\newtheorem{question}{{\bf Question}}
\newenvironment{proof}{{\bf Proof:}}{\hfill$\square$ \medskip}

\begin{document}

\pagestyle{plain}

\title{ The structure of the group of conjugating automorphisms and the linear representation
 of the braid groups of some manifolds
\footnote{Partially supported by the Russian
Foundation for Basic Research (grant number  02--01--01118).}}

\author{Valerij Bardakov}


\maketitle

\begin{abstract}
In this paper we describe the structure of a group of conjugating
automorphisms $C_n$ and prove that this structure is similar to
the structure of a braid group $B_n$ with $n\geq 2$ strings. We
find the linear representation of group $C_n$. Also we prove that
the braid group $B_n(S^2)$ of 2--sphere, mapping class group
$M(0,n)$  of the $n$--punctured 2--sphere and the braid group
$B_3(P^2)$ of the projective plane are linear. Using result of
J.~Dyer, E.~Formanek, E.~Grossman and the faithful linear
representation of Lawrence--Krammer of $B_4$ we construct faithful
linear representation of the automorphism group ${\rm Aut}(F_2)$.
\end{abstract}

\section{Introduction}

\vskip 10pt

It is well known [1,2] that the automorphism group ${\rm Aut}(F_2)$
of a two generated free group $F_2$ can be constructed from cyclic
groups using free and semi direct products. The following problem is open: is it true
for ${\rm Aut}(F_n)$ for $n>2$.
Denote $F'_n=[F_n, F_n]$ (the commutator subgroup of $F_n$) and
$A_n=F_n/F_n'$. Now  $A_n$ is a free abelian group and we have
${\rm Aut}(A_n)\simeq {\rm GL}_n(\mathbb{Z}).$
The natural map from ${\rm Aut}(F_n)$ into ${\rm Aut}(A_n)$ is an
epimorphism. The kernel of this map is called the group
of $IA$--{\it automorphisms} and denoted as ${\rm IA}(F_n)$.
The group ${\rm IA}(F_2)$ is the group of inner automorpisms ${\rm Inn}(F_2)$
which is
isomorphic to $F_2$, but ${\rm Inn}(F_n)$ is a proper subgroup of ${\rm
IA}(F_n)$ if $F_n$ is free of rank greater than 2.

J.~Nielsen showed for $n\leq 3$ and W.~Magnus showed for all $n$ (see [3, chapter 1, \S~4]),
that ${\rm IA}(F_n)$ is generated by the automorphisms
$$
\varepsilon_{ijk} : \left\{
\begin{array}{ll}
x_{i} \longmapsto x_{i}[x_j, x_k] & \mbox{for } k\neq i, j, \\
x_{l} \longmapsto x_{l} & \mbox{for } l\neq i,
\end{array} \right. ~~~~~
\varepsilon_{ij} : \left\{
\begin{array}{ll}
x_{i} \longmapsto x_{j}^{-1}x_ix_j & \mbox{for } i\neq j, \\
x_{l} \longmapsto x_{l} & \mbox{for } l\neq i,
\end{array} \right.
$$
where $[a, b]=a^{-1}b^{-1}ab$ is a commutator $a$ and $b$. J.~Nielsen
showed that ${\rm IA}(F_n)$ is the normal closure in ${\rm
Aut}(F_n)$ of $\varepsilon_{12}.$ The defining relation for ${\rm
IA}(F_n)$, $n\geq 3,$ is not known but S.~Krstic and J.~McCool [4]
proved that ${\rm IA}(F_3)$ is not finite presented.

Denote subgroup of ${\rm IA}(F_n)$ which is generated by $\varepsilon_{ij}$,
$1\leq i\neq j\leq n$, by $Cb_n$. Group $Cb_n$ is called
{\it group of basis--conjugating automorphisms}. J.~McCool [5] proved that this group is
finitely
presented and found defining relations (see below).

Group $Cb_n$ is subgroup of group conjugating automorphisms
$C_n$, where automorphism from ${\rm Aut}(F_n)$
is called {\it conjugating automorphism} if it maps  each generator $x_i$
to $f_i^{-1}x_{\pi (i)}f_i$, where $f_i \in F_n$ and $\pi $ is a permutation
from symmetric group $S_n$. It is clear that if $\pi $ is identical permutation then this
conjugating automorphism lies in $Cb_n$. Automorphisms from $C_n$ which fix product
 $x_1x_2...x_n$ form the braid group $B_n$ on $n$ strings. The braid group $B_n$
has a normal subgroup $P_n$ of finite index. Subgroup $P_n$ is called {\it pure braid group} and quotient $B_n/P_n$
is the symmetric group $S_n$. It is clear that $P_n$ is subgroup of $Cb_n$.
The structure of $P_n$ is well known [6, 7]. It is the semi direct
product of free groups:
$$
P_n=U_n\leftthreetimes (U_{n-1}\leftthreetimes (\ldots \leftthreetimes
(U_3\leftthreetimes U_2))\ldots),
$$
where $U_i\simeq F_{i-1}, i=2,3,\ldots,n.$ The following natural question arises:
is it true  that group of basis--conjugating automorphisms $Cb_n$ is a
semi direct product of some group? In this paper we give
positive answer on this question. It will be proved the following
theorem

\begin{theorem}
The group of basis--conjugating automorphisms $Cb_n$, $n\geq 2,$ is a semi
direct product
$$
Cb_n=D_{n-1}\leftthreetimes (D_{n-2}\leftthreetimes (\ldots
\leftthreetimes (D_2\leftthreetimes D_1))\ldots ),
$$
where $D_i=<\varepsilon_{i+1,1},\varepsilon_{i+1,2},
\ldots,\varepsilon_{i+1,i},
\varepsilon_{1,i+1},\varepsilon_{2,i+1},\ldots,\varepsilon_{i,i+1}>$,
$i=1,2,\ldots,n-1$. Elements $\varepsilon_{i+1,1},\varepsilon_{i+1,2},
\ldots,\varepsilon_{i+1,i}$ generate a free group of rank $i$
and
elements  $\varepsilon_{1,i+1},\varepsilon_{2,i+1},\ldots,\varepsilon_{i,i+1}$
generate a free abelian group of rank $i$.
This decomposition is consistent with the decomposition of $P_n$, i.~e., there is inclusion
$U_{i+1}\leq D_i,$ $i=1,2,\ldots,n-1.$
\end{theorem}

Remark that analogous theorem was proved by  D.~J.~Collins and
N.~D.~Gilbert [8] for group which is the kernel of homomorphism ${\rm Aut}(G) \longrightarrow
{\rm Aut}(\overline{G})$, where $G = G_1 * G_2 * \ldots * G_n$ is
a free product, group $\overline{G} = G_1 \times G_2 \times \ldots \times G_n$
is a direct product and each factor $G_i$ is indecomposible and is not
isomorphic to infinite cyclic group.

As a corollary from this theorem we will find the some normal form
for words which is presented elements from
 $C_n$. We will prove some properties of $C_n$ and $Cb_n$.  We will prove that for $n\geq 2$ the group $C_n$ has no more than 4 generators
and we will find defining relations. For $n\geq 2$ we will prove
that $Cb_n$ has not proper verbal subgroups of finite width. This
result is closely related to the question 14.15 from  "Kourovka
Notebook" [9].

Also one can ask the following question:

\begin{question}
For which $n$  the group generated by $\varepsilon_{ijk}$, $1\leq k\neq i,j\leq n,$
is finitely defined?
\end{question}

One of the most intriguing problems in the theory of braid groups
is the question of whether $B_n$ is linear, i.~e., whether it
admits a faithful representation into a group of matrices over
field. This question was formulated by W.~Burau in 1936. W.~Burau
[10] found an $n$--dimensional linear representation of $B_n$ which
for a long time had been considered as a candidate for faithful
representation. However, as it was established by J.~A.~Moody [11] in
1991 this representation is not faithful for $n\geq 9$. This bound
was improved for $n\geq 6$ by D.~Long and M.~Paton [12] and to
$n\geq 5$ by S.~Bigelow [13]. To this date, it is unknown whether
the Burau representation of $B_4$ is faithful.

In 1990 R.~J.~Lawrence [14] introduced a family of linear representations of $B_n$.
Later D.~Krammer [15] and S.~Bigelow [16] proved that one of this
representation is faithful. Therefore $B_n$ are linear for all $n\geq
2$.

As it was mention above, $B_n$ is subgroup of ${\rm Aut}(F_n).$ But
E.~Formanek and C.~Procesi [17] proved that ${\rm Aut}(F_n)$ is not
linear for all $n\geq 3$. So the following question arises:

\begin{question}
Find the maximal related inclusion subgroup of ${\rm Aut}(F_n)$, $n\geq 3$,
which contains  $B_n$ and is linear. In particular, is it true that
group $C_n$ is linear?
\end{question}

The second part of this problem was formulated in [9, Problem
15.9]. The one natural way for answer on this question is
to extend the known representations of $B_n$ to a representation of $C_n$.

In this paper we define an extension of the Burau representation on
$C_n$ and prove that the linear representation of Lawrence--Krammer
extending on $C_3$ and for $n\geq 4$ we define the extension of
this representation for some restrictions on parameters of
representation.

Braid group $B_n$ is a particular case of general construction braid
group $B_n(M)$ of manifold $M$ on $n$ strings.

\begin{question}
For which manifold $M$ and natural $n$ the braid group $B_n(M)$
on $n$ strings is linear?
\end{question}

The greatest interest in this problem is connected to manifolds of
dimension two. In the braid groups of these manifolds only $B_n(S^2)$
of 2--sphere $S^2$ and $B_n(P^2)$
of projective plane $P^2$ have a torsion.

In this article we will prove that the mapping class group $M(0,n)$ of the
$n$--punctured 2--sphere and braid group $B_n(S^2)$
of 2--sphere are linear for all $n\geq 2$. Also we will prove that braid group
 $B_3(P^2)$ of projective plane is linear (for $n=1, 2$ this group is finite).
In all cases we will construct corresponding linear
representation. The result that $B_n(S^2)$ is linear was announced
in [18]. As it was known to the author S.~Bigelow and R.~D.~Budney
[19] give another proof that groups $B_n(S^2)$ and $M(0,n)$ are
linear.

As we noted above group ${\rm Aut}(F_n)$ is not linear for $n\geq
3$. It was proved in [20] that ${\rm Aut}(F_2)$ is linear if and
only if braid group $B_4$ on 4 strings is linear.
In last section we will construct the faithful linear representation of ${\rm Aut}(F_2)$.

I thankful to J. Meier for pointing out on work [8] and on some
non-corrected statement in the first version of the manuscript.

\vskip 20pt


\section{Definitions and notations}

\vskip 12pt

We remind some known facts which can be found in [21, 6, 7].

The group $G$ is a {\it semi direct product} of $A$ and $B$ if there exist subgroups
$H$ and $K$ in $G$ such that
$$
G=HK,~~~A\simeq H\unlhd G,~~~B\simeq K,~~~H\cap K=1.
$$
Semi direct product is denoted as $G=A\leftthreetimes B$. If $A$ and
$B$ have presentations
$$A=<a_1, a_2,\ldots,a_k \parallel A_1,
A_2,\ldots,A_p>,~~~
B=<b_1, b_2,\ldots,b_l \parallel B_1, B_2,\ldots,B_q>,
$$
then $G=A\leftthreetimes B$ have presentation
$$
G=<a_1, a_2,\ldots,a_k; b_1, b_2,\ldots,b_l \parallel A_1, A_2,\ldots,A_p;
B_1, B_2,\ldots,B_q; b_i^{-1}a_jb_i=C_{ij},
$$
$$
~~~~~~~~~~~~~~~~~~~~~~~~~~~~~~~~~~~~~~~~~~~~~~~~~~1\leq i\leq l,~~~ 1\leq j\leq k>,
$$
where $A_1, A_2,\ldots,A_p, C_{ij}$ are words in alphabet
$\{ a_1^{\pm 1},a_2^{\pm 1},\ldots,a_k^{\pm 1}\}$, $B_1, B_2,\ldots,B_q$
are words in alphabet
$\{ b_1^{\pm 1},b_2^{\pm 1},\ldots,b_l^{\pm 1}\},$ elements $C_{i1}, C_{i2},\ldots,C_{ik}$
generate group $A$ for each $i\in \{ 1,2,\ldots,l \}$ and conjugation by element $b_i$ induces
automorphism of group $A$.

Let $M$ be a manifold of dimension $\geq 2$. {\it Configuration space} $F_n(M)$,
$n\in \mathbb{N}$, for manifold $M$ is a set
$$
F_n(M)=\{ (z_1, z_2,\ldots,z_n)\in M^n \vert z_i\neq z_j~~
\mbox{for}~~i\neq j \}
$$
of ordering collections of $n$ distinct points from $M$. The fundamental group
$P_n(M)=\pi_1(F_n(M))$ of the space $F_n(M)$ is the {\it pure braid group} with $n$
strings of the manifold $M$. Symmetric group $S_n$ acts on the space $F_n(M)$ permute the
coordinats. This action is free and induce the regular covering of orbit space $F_n(M)/S_n$
by $F_n(M)$. The fundamental group $B_n(M)=\pi_1(F_n(M)/S_n)$ is called the {\it full
braid group} of $M$, or more simply the {\it braid group} of $M$.
The regular covering projection $F_n(M)\longrightarrow F_n(M)/S_n$
is induce the short exact sequence
$$
1\longrightarrow P_n(M)\longrightarrow B_n(M)\longrightarrow
S_n\longrightarrow 1.
$$
If $M$ is a closed, smooth manifold of dimension $n\geq 2$, then the inclusion map
 $F_n(M)\longrightarrow M^n$
induces a surjective homomorphism $P_n(M)\longrightarrow \pi_1(M)\times \ldots \times \pi_1(M)$
pure braid group $P_n(M)$ on direct product $n$ copies of fundamental group $\pi_1(M)$.
If $\mbox{dim}M>2$ then this homomorphism also injective. The braid groups of
manifolds of dimension 2 represent
the largest interest.

Among the 2--dimension manifolds special role play the 2--sphere $S^2$
and the projective plane $P^2$ because the braid groups only this manifolds have elements of finite order.
If $M$ is a closed surface different from $S^2$ or $P^2$, then in
the following sequence of groups
$$
1\longrightarrow P_n(E^2)\longrightarrow P_n(M)\longrightarrow
\prod^n_{i=1}\pi_1(M)\longrightarrow 1
$$
the kernel of each homomorphism is equal to the normal
closure of the image of the previous homomorphism in the sequence.
In this sequence $E^2$ denotes the Euclidean plane and $\prod^n_{i=1}\pi_1(M)$
denoted the direct product of $n$ copies of $\pi_1(M)$.

The classical braid group of Artin $B_n$ is the braid group $B_n(E^2)$ of Euclidean
plane $E^2$. We will call $B_n$  simply braid group.

The braid group $B_n$, $n\geq 2$, with $n$ strings can be
defined as the group generated by $n-1$ generators
$\sigma_1,\sigma_2,...,\sigma_{n-1}$ with defining relations
$$
\sigma_i\sigma_{i+1}\sigma_i=\sigma_{i+1}\sigma_i\sigma_{i+1},~~~ i=1,2,...,n-2, \eqno{(1)}
$$
$$
\sigma_i \sigma_j = \sigma_j \sigma_i,~~~ |i-j|\geq 2. \eqno{(2)}
$$
There is a homomorphism $\nu : B_n\longrightarrow S_n$ from the
group $B_n$ to the symmetric group $S_n$ on $n$ letters defined by
$$
\nu(\sigma_i)=(i,i+1),~~~i=1,2,\ldots,n-1.
$$
The kernel of homomorphism $\nu $ is the pure braid group $P_n$.
The group $P_n$ admits a presentation with generators
$$
a_{i,i+1}=\sigma_i^2,
$$
$$
a_{ij}=\sigma_{j-1}\sigma_{j-2}\ldots\sigma_{i+1}\sigma_i^2\sigma_{i+1}^{-1}\ldots
\sigma_{j-2}^{-1}\sigma_{j-1}^{-1},~~~i+1< j \leq n.
$$
and defining relations:
$$
\begin{array}{ll}
a_{ik}^{-\nu }a_{kj}a_{ik}^{\nu }=\left( a_{ij}a_{kj}\right) ^{\nu }a_{kj}
\left( a_{ij}a_{kj}\right) ^{-\nu }, & \\
& \\
a_{km}^{-\nu }a_{kj}a_{km}^{\nu }=\left( a_{kj}a_{mj}\right) ^{\nu }a_{kj}
\left( a_{kj}a_{mj}\right) ^{-\nu }, &~~~m < j,\\
& \\
a_{im}^{-\nu }a_{kj}a_{im}^{\nu }=\left[ a_{ij}^{-\nu }, a_{mj}^{-\nu }\right] ^{\nu }a_{kj}
\left[ a_{ij}^{-\nu }, a_{mj}^{-\nu }\right] ^{-\nu }, & ~~~i < k <
m,\\
& \\
a_{im}^{-\nu }a_{kj}a_{im}^{\nu }=a_{kj}, & ~~~k < i;~~~m < j~~
\mbox{or}~~ m < k,
\end{array}
$$
where $\nu=\pm 1$.

The subgroup $U_n$ of $P_n$ which is generated by $a_{1,n},
a_{2,n},\ldots,a_{n-1,n}$ is free and normal in $P_n$. The group $P_n$ is a
semi direct product of $U_n$ and $P_{n-1}$. Hence the group $P_n$
is a semi direct product
$$
P_n=U_n\leftthreetimes (U_{n-1}\leftthreetimes (\ldots \leftthreetimes
(U_3\leftthreetimes U_2))\ldots),
$$
where $U_i=<a_{1,i},
a_{2,i},\ldots,a_{i-1,i}>, i=2,3,\ldots,n,$ is a free group of
rank $i-1$.

The braid group $B_n$ has a faithful representation as a group of
automorphisms of a free group $F_n=<x_1, x_2,\ldots, x_n>$ of rank
$n$. The representation is induced by a mapping from $B_n$ to ${\rm Aut}(F_n)$
defined by
$$
\sigma_{i} : \left\{
\begin{array}{ll}
x_{i} \longmapsto x_{i}x_{i+1}x_i^{-1}, &  \\ x_{i+1} \longmapsto
x_{i}, & \\ x_{l} \longmapsto x_{l} & ~~ l\neq i,i+1.
\end{array} \right.
$$
The generator $a_{rs}$ of $P_n$ defines the following automorphism
$$
a_{rs} : \left\{
\begin{array}{ll}
x_{i} \longmapsto x_{i} & \mbox{if }~~ s < i~~ \mbox{or}~~ i < r, \\
& \\
x_{r} \longmapsto x_{r}x_{s}x_rx_s^{-1}x_r^{-1}, &  \\
& \\
x_{i} \longmapsto [x_{r}^{-1}, x_s^{-1}]x_i[x_{r}^{-1}, x_s^{-1}]^{-1} & \mbox{if }~~ r < i <s,\\
& \\
x_{s} \longmapsto x_{r}x_sx_r^{-1}. &
\end{array} \right.
$$

By the theorem of Artin [7, Theorem 1.9] automorphism $\beta \in B_n \subset {\rm Aut}(F_n)$
if and only if $\beta $ satisfies to the two conditions
$$
1)~~~ \beta(x_i)=f_i^{-1}x_{\pi(i)}f_i,~~~1\leq i\leq n,
$$
$$
2)~~~ \beta(x_1x_2 \ldots x_n)=x_1x_2 \ldots x_n,
$$
where $\pi $ is a permutation of  $\{1, 2, \ldots, n \}$ and $f_i=f_i(x_1,x_2,\ldots,x_n)$
is a word in the generators of $F_n$.

The automorphism of $F_n$ calls {\it conjugating automorphism} if it
satisfies condition 1). Let $C_n$ be a group of conjugating
automorphisms. Obvious, that it has the normal subgroup $Cb_n$ with
factor group $C_n/Cb_n$ isomorphic to the symmetric group $S_n$. The
group $Cb_n$ calls the group of {\it basis--conjugating
automorphisms}. J.~McCool [5] proved, that the group $Cb_n$ is
generated by automorphisms
$$
\varepsilon_{ij} : \left\{
\begin{array}{ll}
x_{i} \longmapsto x_{j}^{-1}x_ix_j &  i\neq j, \\
x_{l} \longmapsto x_{l} &  l\neq i,
\end{array} \right.
$$
$i\leq i\neq j \leq n$, with defining relations
$$
\begin{array}{lr}
\varepsilon_{ij}\varepsilon_{kl}=\varepsilon_{kl}\varepsilon_{ij},
&~~~~~~~~~~~~~~~~~~~~~~~~~~~~~~~~~~~~~~~~~~~~~~~~~~~~~~~~  (3) \\
& \\
\varepsilon_{ij}\varepsilon_{kj}=\varepsilon_{kj}\varepsilon_{ij},
&~~~~~~~~~~~~~~~~~~~~~~~~~~~~~~~~~~~~~~~~~~~~~~~~~~~~~~~~ (4) \\
& \\
(\varepsilon_{ij}\varepsilon_{kj})\varepsilon_{ik}=\varepsilon_{ik}(\varepsilon_{ij}
\varepsilon_{kj}).
&~~~~~~~~~~~~~~~~~~~~~~~~~~~~~~~~~~~~~~~~~~~~~~~~~~~~~~~~ (5) \\
\end{array}
$$
where distinct letters denote distinct indexes.

\vskip 20pt

\section{Decomposition of the group of basis--conjugating automorphisms}

\vskip 12pt

In this section it will be proved that the structure of $Cb_n$ is similar to the
structure of pure braid group $P_n$. The result is the following:

\begin{theorem}
The group of basis--conjugating automorphisms $Cb_n$, $n\geq 2$ is a semi
direct product
$$
Cb_n=D_{n-1}\leftthreetimes (D_{n-2}\leftthreetimes (\ldots
\leftthreetimes (D_2\leftthreetimes D_1))\ldots ),
$$
where $D_i=<\varepsilon_{i+1,1},\varepsilon_{i+1,2},
\ldots,\varepsilon_{i+1,i},
\varepsilon_{1,i+1},\varepsilon_{2,i+1},\ldots,\varepsilon_{i,i+1}>$,
$i=1,2,\ldots,n-1$. Elements $\varepsilon_{i+1,1},\varepsilon_{i+1,2},
\ldots,\varepsilon_{i+1,i}$ generate a free group of rank $i$
and
elements  $\varepsilon_{1,i+1},\varepsilon_{2,i+1},\ldots,\varepsilon_{i,i+1}$
generate a free abelian group of rank $i$.
This decomposition is consistent with the decomposition of $P_n$, i.~e., there is inclusion
$U_{i+1}\leq D_i,$ $i=1,2,\ldots,n-1.$
\end{theorem}

As an immediate consequence of this theorem we obtain the normal
form of words in $C_n$.

\begin{corollary}
Every element $\beta \in C_n$ can be written  in the form
$$
\beta=\beta_1\beta_2\ldots \beta_{n-1}\pi_{\beta},
$$
where $\pi_{\beta}$ is a permutation automorphism and each
$\beta_j$ belongs to the subgroup $D_j$.
\end{corollary}

For prove of the theorem 3.1 we need the following lemma which
follow from defining relations of $Cb_n$

\begin{lemma}
The following conjugating rules are true in $Cb_n$:
$$
\begin{array}{ll}
1)~~& \varepsilon_{ij}^{-\nu }\varepsilon_{kl}\varepsilon_{ij}^{\nu }=\varepsilon_{kl}, \\
\\
2)~~& \varepsilon_{ij}^{-\nu }\varepsilon_{kj}\varepsilon_{ij}^{\nu }=\varepsilon_{kj},    \\
\\
3)~~& \varepsilon_{ij}^{-\nu }\varepsilon_{ki}\varepsilon_{ij}^{\nu }=
\varepsilon_{kj}^{\nu }\varepsilon_{ki}\varepsilon_{kj}^{-\nu },  \\
\\
4)~~& \varepsilon_{ij}^{-\nu }\varepsilon_{ik}\varepsilon_{ij}^{\nu }=
\varepsilon_{kj}^{\nu }\varepsilon_{ik}\varepsilon_{kj}^{-\nu },  \\
\\
5)~~& \varepsilon_{ij}^{-\nu }\varepsilon_{jk}\varepsilon_{ij}^{\nu }=
[\varepsilon_{kj}^{-\nu }, \varepsilon_{ik}]\varepsilon_{jk}.  \\
\end{array}
$$
\end{lemma}
where $\nu =\pm 1$ and distinct letters denote distinct indexes.

Define the following subgroups in $Cb_n$:
$$
D_i=<\varepsilon_{i+1,1}, \varepsilon_{i+1,2},\ldots,\varepsilon_{i+1,i},
\varepsilon_{1,i+1},
\varepsilon_{2,i+1},\ldots,\varepsilon_{i,i+1}>,~~~i=1,2,\ldots,n-1.
$$

\begin{lemma} Let $k$ and $l$ be the integer numbers such that $2 \leq k < l \leq n$.
Then the group $D_{k-1}$ lies in normalizer of
$D_{l-1}$.
\end{lemma}

\begin{proof}
From definition we have
$$
D_{k-1}=<\varepsilon_{k,1}, \varepsilon_{k,2},\ldots,\varepsilon_{k,k-1},
\varepsilon_{1,k},
\varepsilon_{2,k},\ldots,\varepsilon_{k-1,k}>,
$$
$$
D_{l-1}=<\varepsilon_{l,1}, \varepsilon_{l,2},\ldots,\varepsilon_{l,l-1},
\varepsilon_{1,l},
\varepsilon_{2,l},\ldots,\varepsilon_{l-1,l}>.
$$
Let $\varepsilon_{k,i}$, $1\leq i \leq k-1$ be some generator of
$D_{k-1}$.
We should prove that element $\varepsilon^{-\nu }_{k,i}d\varepsilon^{-\nu }_{k,i},$
$\nu =\pm 1$, lies in $D_{l-1}$ for each  generator $d$ of
$D_{l-1}$. By Lemma 1 we have the following conjugating rules:
$$
\begin{array}{l}
\varepsilon_{ki}^{-\nu }\varepsilon_{lj}\varepsilon_{ki}^{\nu }=\varepsilon_{lj},~~~
\varepsilon_{ki}^{-\nu }\varepsilon_{jl}\varepsilon_{ki}^{\nu }=\varepsilon_{jl},~~~j\neq i,
k,
 \\
 \\
 \varepsilon_{ki}^{-\nu }\varepsilon_{li}\varepsilon_{ki}^{\nu }=\varepsilon_{li},    \\
\\
 \varepsilon_{ki}^{-\nu }\varepsilon_{lk}\varepsilon_{ki}^{\nu }=
\varepsilon_{li}^{\nu }\varepsilon_{lk}\varepsilon_{li}^{-\nu },  \\
\\
 \varepsilon_{ki}^{-\nu }\varepsilon_{kl}\varepsilon_{ki}^{\nu }=
\varepsilon_{li}^{\nu }\varepsilon_{kl}\varepsilon_{li}^{-\nu },  \\
\\
 \varepsilon_{ki}^{-\nu }\varepsilon_{il}\varepsilon_{ki}^{\nu }=
[\varepsilon_{li}^{-\nu }, \varepsilon_{kl}]\varepsilon_{il}.  \\
\end{array}
$$
We see that right parts of this rules lie in $D_{l-1}$. Hence in this case lemma is true.

Let $\varepsilon_{ik}$, $1\leq i \leq k-1$ be a generator of $D_{k-1}$.
By Lemma 1 we have the following conjugating rules:
$$
\begin{array}{l}
\varepsilon_{ik}^{-\nu }\varepsilon_{lj}\varepsilon_{ik}^{\nu }=\varepsilon_{lj},~~~
\varepsilon_{ik}^{-\nu }\varepsilon_{jl}\varepsilon_{ik}^{\nu }=\varepsilon_{jl},~~~j\neq i,
k,
 \\
 \\
 \varepsilon_{ik}^{-\nu }\varepsilon_{lk}\varepsilon_{ik}^{\nu }=\varepsilon_{lk},    \\
\\
 \varepsilon_{ik}^{-\nu }\varepsilon_{li}\varepsilon_{ik}^{\nu }=
\varepsilon_{lk}^{\nu }\varepsilon_{li}\varepsilon_{lk}^{-\nu },  \\
\\
 \varepsilon_{ik}^{-\nu }\varepsilon_{il}\varepsilon_{ik}^{\nu }=
\varepsilon_{lk}^{\nu }\varepsilon_{il}\varepsilon_{lk}^{-\nu },  \\
\\
 \varepsilon_{ik}^{-\nu }\varepsilon_{kl}\varepsilon_{ik}^{\nu }=
[\varepsilon_{lk}^{-\nu }, \varepsilon_{il}]\varepsilon_{kl}.  \\
\end{array}
$$
We see again that right parts of this rules lie in $D_{l-1}$.
Hence the lemma is true in all cases.
\end{proof}

From this lemma with the help of induction we have

\begin{corollary}
Subgroup $D_{n-1}$ is normal in $Cb_n$.
\end{corollary}

The following lemma is a part of Theorem 3.1.

\begin{lemma} If $n\geq 2$ then the group of basis--conjugating
automorphisms $Cb_n$ is a semi direct product
$$
Cb_n=D_{n-1}\leftthreetimes (D_{n-2}\leftthreetimes \ldots (\leftthreetimes (D_{2}\leftthreetimes
D_{1}))\ldots),
$$
where $D_i=<\varepsilon_{i+1,1}, \varepsilon_{i+1,2},\ldots,\varepsilon_{i+1,i},
\varepsilon_{1,i+1},
\varepsilon_{2,i+1},\ldots,\varepsilon_{i,i+1}>,$ $
<\varepsilon_{i+1,1}, \varepsilon_{i+1,2},\ldots,\varepsilon_{i+1,i}>\simeq F_i,$
$<\varepsilon_{1,i+1},
\varepsilon_{2,i+1},\ldots,\varepsilon_{i,i+1}>\simeq \mathbb{Z}^i,$ $i=1,2,\ldots,n-1.$
\end{lemma}

\begin{proof}
For $n=2$ we have
$$
Cb_2=D_1=<\varepsilon_{21}, \varepsilon_{12}>.
$$
Since the group $Cb_2$ has not relations, then $D_1\simeq F_2=F_1*\mathbb{Z}$
and we have the basis of induction.

We will prove that $Cb_n=D_{n-1}\leftthreetimes Cb_{n-1}$. From Corollary 3.2
subgroup
$D_{n-1}$ is normal in $Cb_n$. We divide  defining relations (3)--(5) of $Cb_n$
on three parts which is not intersecting.  In the first part we pose defining relations which contain
only generators from $D_{n-1}$. We see that in every pair of generators of $D_{n-1}$
one index is equal to $n$. Hence in the first part of defining relations  only relations from
 (4) lie, i. e.
$$
\varepsilon_{in}\varepsilon_{jn}=\varepsilon_{jn}\varepsilon_{in},~~~1 \leq i\neq
j\leq n-1. \eqno{(6)}
$$
Therefore  $\varepsilon_{1n}, \varepsilon_{2n},\ldots,\varepsilon_{n-1,n}$
generate the abelian subgroup of $D_{n-1}$.
Also, it is not hard to check that elements $\varepsilon_{n1}, \varepsilon_{n2},\ldots,\varepsilon_{n,n-1}$
are the free generators of free group $F_{n-1}$.
 In the second part
of relations  relations lie only in generators of $Cb_{n-1}$. The third part
of relations relations  in generators of $D_{n-1}$ and $Cb_{n-1}$ simultaneously lie.
We should prove that this relations define an action of $Cb_{n-1}$ on $D_{n-1}$.

Consider relations from (3) which include generators of $D_{n-1}$ and $Cb_{n-1}$ simultaneously.
This relations look as
$$
\varepsilon_{ij}\varepsilon_{kn}=\varepsilon_{kn}\varepsilon_{ij},~~~
\varepsilon_{ij}\varepsilon_{nk}=\varepsilon_{nk}\varepsilon_{ij},~~~1 \leq
i,j,k < n. \eqno{(7)}
$$
From these relations we have the following conjugating rules
$$
\varepsilon_{ij}^{-\nu }\varepsilon_{kn}\varepsilon_{ij}^{\nu }=\varepsilon_{kn},~~~
\varepsilon_{ij}^{-\nu }\varepsilon_{nk}\varepsilon_{ij}^{\nu }=\varepsilon_{nk},~~~
\nu =\pm 1. \eqno{(8)}
$$

If defining relation from (4) includes generators of $D_{n-1}$
and $Cb_{n-1}$, then it look as
$$
\varepsilon_{nj}\varepsilon_{kj}=\varepsilon_{kj}\varepsilon_{nj},~~~1 \leq k,j < n.
$$
From this relations we have
$$
\varepsilon_{kj}^{-\nu }\varepsilon_{nj}\varepsilon_{kj}^{\nu }=\varepsilon_{nj},
~~~1 \leq k,j < n~~~
\nu =\pm 1. \eqno{(9)}
$$

If relation from (5) includes generators of $D_{n-1}$
and $Cb_{n-1}$ simultaneously, then this relation has on of the
following three types
$$
\begin{array}{l}
(\varepsilon_{ij}\varepsilon_{nj})\varepsilon_{in}=
\varepsilon_{in}(\varepsilon_{ij}\varepsilon_{nj}), \\
\\
(\varepsilon_{nj}\varepsilon_{ij})\varepsilon_{ni}=
\varepsilon_{ni}(\varepsilon_{nj}\varepsilon_{ij}),  \\
\\
(\varepsilon_{in}\varepsilon_{jn})\varepsilon_{ij}=
\varepsilon_{ij}(\varepsilon_{in}\varepsilon_{jn}),~~~1\leq i,j < n.  \\
\end{array}
$$
From this relations we have the formulas of conjugating:
$$
\begin{array}{lr}
\varepsilon_{ij}^{-\nu }\varepsilon_{in}\varepsilon_{ij}^{\nu }=
\varepsilon_{nj}^{\nu }\varepsilon_{in}\varepsilon_{nj}^{-\nu },
&~~~~~~~~~~~~~~~~~~~~~~~~~~(10) \\
& \\
\varepsilon_{ij}^{-\nu }\varepsilon_{ni}\varepsilon_{ij}^{\nu }=
\varepsilon_{nj}^{\nu }\varepsilon_{ni}\varepsilon_{nj}^{-\nu },
&~~~~~~~~~~~~~~~~~~~~~~~~~~(11) \\
& \\
\varepsilon_{ij}^{-\nu }\varepsilon_{jn}\varepsilon_{ij}^{\nu }=
[\varepsilon_{nj}^{-\nu },
\varepsilon_{in}]\varepsilon_{jn},~~~1\leq i,j < n,~~~\nu =\pm 1.
&~~~~~~~~~~~~~~~~~~~~~~~~~~~~~~~~~~~~~~~~~~~~~~~~~~(12) \\
\end{array}
$$
So, defining relations of $Cb_n$ which include generators of $D_{n-1}$
and $Cb_{n-1}$ simultaneously are equivalent to (8)--(12). Hence,
we define action of
$Cb_{n-1}$ on $D_{n-1}$. Since
$D_{n-1}\bigcap Cb_{n-1}=1$ then $Cb_n=D_{n-1}\leftthreetimes Cb_{n-1}$.
From inductive hypothesis we have the needed decomposition.
\end{proof}

The group of conjugating automorphisms $C_n$ include the group of
basis--conjugating automorphisms $Cb_n$. A.~G.~Savushkina [22]
proved that  $C_n$ is a semi direct product $C_n=Cb_n\leftthreetimes S_n$ of $Cb_n$
and symmetric group  $S_n$, where $S_n$  generating by
automorphisms:
$$
\alpha_{i} : \left\{
\begin{array}{ll}
x_{i} \longmapsto x_{i+1}, &  \\
x_{i+1} \longmapsto x_{i}, & \\
x_{l} \longmapsto x_{l}, & l\neq i,i+1,\\
\end{array} \right.
$$
$i=1,2,\ldots, n-1,$ and defining relations:
$$
\alpha_j^2=1,~~~j=1,2,\ldots, n-1, \eqno{(13)}
$$
$$
\alpha_j\alpha_{j+1}\alpha_j=\alpha_{j+1}\alpha_{j}\alpha_{j+1},~~~j=1,2,\ldots,n-2, \eqno{(14)}
$$
$$
\alpha_k\alpha_{l}=\alpha_{l}\alpha_{k},~~~1\leq k,l\leq n-1,~~~|k-l|\geq 2. \eqno{(15)}
$$
Whole group $C_n$ is generated by $\varepsilon_{ij}$, $1\leq i\neq j\leq n$;
$\alpha_{k}$, $1\leq k\leq n-1$ and defined by relations  (3)--(5), (13)--(15)
and
$$
\varepsilon_{ij}\alpha_{k}=\alpha_{k}\varepsilon_{ij},~~~k\neq
i-1,i,j-1,j,  \eqno{(16)}
$$
$$
\varepsilon_{ij}\alpha_{i}=\alpha_{i}\varepsilon_{i+1,j},~~~
\varepsilon_{ij}\alpha_{j}=\alpha_{j}\varepsilon_{i,i+1},~~~
\varepsilon_{i,i+1}\alpha_{i}=\alpha_{i}\varepsilon_{i+1,i},~~~
  \eqno{(17)}
$$
(see [22, Lemma 1]).

Evidently, pure braid group $P_n$ is subgroup of $Cb_n$.
We can express generators of $P_n$ over generators of $Cb_n$.

\begin{lemma} Generators of $P_n$ can be express as
$$
a_{i,i+1}=\varepsilon_{i,i+1}^{-1}\varepsilon_{i+1,i}^{-1},~~~i=1,2,\ldots,n-1, \eqno{(18)}
$$
$$
a_{ij}=\varepsilon_{j-1,i}\varepsilon_{j-2,i}\ldots\varepsilon_{i+1,i}
(\varepsilon_{ij}^{-1}\varepsilon_{ji}^{-1})\varepsilon_{i+1,i}^{-1}\ldots
\varepsilon_{j-2,i}^{-1}\varepsilon_{j-1,i}^{-1}=
$$
$$
=\varepsilon_{j-1,j}^{-1}\varepsilon_{j-2,j}^{-1}\ldots\varepsilon_{i+1,j}^{-1}
(\varepsilon_{ij}^{-1}\varepsilon_{ji}^{-1})\varepsilon_{i+1,j}\ldots
\varepsilon_{j-2,j}\varepsilon_{j-1,j},~~~2 \leq i+1 < j \leq n. \eqno{(19)}
$$
\end{lemma}

\begin{proof} The formula (18) follows from actions of $P_n$ and $Cb_n$ on
free group $F_n$.
For proof of the first part of (19) we should verify that
automorphism
$$
a_{ij}\varepsilon_{ji}
(\varepsilon_{j-1,i}\varepsilon_{j-2,i}\ldots\varepsilon_{i+1,i})
\varepsilon_{ij}(\varepsilon_{i+1,i}^{-1}\ldots\varepsilon_{j-2,i}^{-1}\varepsilon_{j-1,i}^{-1})
$$
is identical. For proof of the second part of (19) we should verify
relation
$$
\varepsilon_{ki}(\varepsilon_{ij}^{-1}\varepsilon_{ji}^{-1})\varepsilon_{ki}^{-1}=
\varepsilon_{kj}^{-1}(\varepsilon_{ij}^{-1}\varepsilon_{ji}^{-1})\varepsilon_{kj},~~~
i+1\leq k\leq j-1. \eqno{(20)}
$$
From Lemma 3.1 we have
$$
\varepsilon_{ki}\varepsilon_{ij}^{-1}\varepsilon_{ki}^{-1}=\varepsilon_{ij}^{-1}
[\varepsilon_{kj}, \varepsilon_{ji}],~~~
\varepsilon_{ki}\varepsilon_{ji}^{-1}\varepsilon_{ki}^{-1}=\varepsilon_{ji}^{-1}.
$$
We use this formulas and rewrite the left--hand side of (20) in
the following form
$$
\varepsilon_{ki}(\varepsilon_{ij}^{-1}\varepsilon_{ji}^{-1})\varepsilon_{ki}^{-1}=
\varepsilon_{ij}^{-1}
[\varepsilon_{kj}, \varepsilon_{ji}]\varepsilon_{ji}^{-1}=
\varepsilon_{ij}^{-1}\varepsilon_{kj}^{-1}\varepsilon_{ji}^{-1}\varepsilon_{kj}.
$$
Since $\varepsilon_{ij}^{-1}$ and
$\varepsilon_{kj}^{-1}$ are commutative, then we have (20). From this relation we have required
formula.
\end{proof}

Now we can complete the proof of the Theorem 3.1. We know that
 free group $U_{i+1}$ generating by
$a_{1,i+1}, a_{2,i+1},\ldots,a_{i,i+1}$.
By Lemma 3.4 we have
$$
a_{k,i+1}=
\varepsilon_{i,i+1}^{-1}\varepsilon_{i-1,i+1}^{-1}\ldots\varepsilon_{k+1,i+1}^{-1}
(\varepsilon_{k,i+1}^{-1}\varepsilon_{i+1,k}^{-1})\varepsilon_{k+1,i+1}\ldots
\varepsilon_{i-1,i+1}\varepsilon_{i,i+1},~~~1 \leq k \leq i-1,
$$
$$
a_{i,i+1}=\varepsilon_{i,i+1}^{-1}\varepsilon_{i+1,i}^{-1}.
$$
Hence $U_{i+1}$ is subgroup of $D_i$.

If we use the Tit's transformation, then not hard to prove that $D_i$ generating by
 $a_{1,i+1}, a_{2,i+1},\ldots,a_{i,i+1};$
$\varepsilon_{1,i+1},
\varepsilon_{2,i+1},\ldots,\varepsilon_{i,i+1}.$
The proof of the theorem is complete.

This theorem reduces the description of $Cb_n$ to the description
of its subgroups $D_{n-1}, \ldots, D_1.$ We can see that this
subgroups (with the exception of $D_1$) have the complicated
structure. In [8, p. 167] the hypothesis was formulated that this
subgroups are not finitely define.

By Theorem 3.1 we have $Cb_3 = D_2 \leftthreetimes D_1,$ where $D_2 = < \varepsilon_{31},
\varepsilon_{32}, \varepsilon_{13}, \varepsilon_{23}>$, $D_1 = < \varepsilon_{21},
\varepsilon_{12}> \simeq F_2.$ Elements $\varepsilon_{13}$ and
$\varepsilon_{23}$ from $D_2$ are commute. Let us prove that there
are many another relations in $D_2$. As it was noted in proof of
Theorem 3.1 conjugation by element from $D_1$ induses automorphism
of $D_2.$ This implies that if we conjugate relation $\varepsilon_{13} \varepsilon_{23} =
\varepsilon_{23} \varepsilon_{13}$ then we get relation from
$D_2.$ In particular, conjugating relation $\varepsilon_{13} \varepsilon_{23} =
\varepsilon_{23} \varepsilon_{13}$ by different powers of $\varepsilon_{12}$ and $\varepsilon_{21}$
we get the following relations
$$
[\varepsilon_{13} \varepsilon_{23}, \varepsilon_{32}^k
\varepsilon_{13}  \varepsilon_{32}^{-k}] = 1,
$$
$$
[\varepsilon_{13} \varepsilon_{23}, \varepsilon_{31}^k
\varepsilon_{23}  \varepsilon_{31}^{-k}] = 1,
$$
for each integer $k$.
\vskip 20pt

\section{Some property of  group of conjugating automorphisms}

\vskip 12pt

We find the abelianized groups $C_n/C_n'$ and $Cb_n/Cb_n'$ from the following known fact.

\begin{lemma} ([22, Theorem 1]). Group $C_n$ is generated by automorphisms $\sigma_i$,
$\alpha_i$, $i=1,2,\ldots,n-1$ and determined by relations (1)--(2), (13)--(15)
and
$$
\alpha_i\sigma_j=\sigma_j\alpha_i,~~~~~~~|i-j|\geq 2, \eqno{(21)}
$$
$$
\sigma_i\alpha_{i+1}\alpha_i=\alpha_{i+1}\alpha_i\sigma_{i+1},
~~~\sigma_i\sigma_{i+1}\alpha_i=\alpha_{i+1}\sigma_i\sigma_{i+1}. \eqno{(22)}
$$
\end{lemma}

\begin{proposition} a) The abelianized group $C_n/C_n'$, $n\geq 2,$ is isomorphic to
 $\mathbb{Z}_2\times \mathbb{Z}$.

b) The abelianized group $Cb_n/Cb_n'$, $n\geq 2,$ is isomorphic to $\mathbb{Z}^{n(n-1)}$.
\end{proposition}

The group of inner automorphisms ${\rm Inn}(F_n)$ of $F_n$ is subgroup of $Cb_n$.

\begin{proposition} Intersection $B_n\bigcap {\rm Inn}(F_n)$, $n\geq 2$, is
infinite cyclic group which generating by inner automorphisms that is conjugation
by $x_1x_2\ldots x_n$.
\end{proposition}

\begin{proof} immediately follows from Artin's theorem (see \S~2) which
gives description of $B_n$ as subgroup of ${\rm Aut}(F_n)$. Elements from $B_n$
fixing the product $x_1x_2\ldots x_n$. Inner automorphisms which has this
property is only conjugation by $(x_1x_2\ldots x_n)^k$, $k\in \mathbb{Z}$. This
completes the proof.
\end{proof}

We will show that $C_n$, $n\geq 2$, can be generated at most 4 elements. By the Lemma
4.1
we have
$$
C_2=<\alpha_1, \sigma_1 \parallel \alpha_1^2=1>\simeq
\mathbb{Z}_2 * \mathbb{Z}.
$$
For $n\geq 3$ we have

\begin{proposition} The group $C_n$, $n\geq 3$, can be generated by
 $\alpha_1$, $\alpha $, $\sigma_1$, $\sigma $ with defining
 relations
$$
\sigma^n=(\sigma \sigma_1)^{n-1},~~~\sigma_1(\sigma^{-j}\sigma_1\sigma^j)=
(\sigma^{-j}\sigma_1\sigma^j)\sigma_1,~~~2\leq j\leq n/2, \eqno{(23)}
$$
$$
\alpha^2=1,~~~\alpha^n=(\alpha_1\alpha )^{n-1},~~~\alpha_1(\alpha^{-j}\alpha_1\alpha^j)=
(\alpha^{-j}\alpha_1\alpha^j)\alpha_1,~~~2\leq j\leq n/2, \eqno{(24)}
$$
$$
\alpha^{-i}\alpha_1\alpha^i\sigma^{j}\sigma_1\sigma^{-j}=
\sigma^{j}\sigma_1\sigma^{-j}\alpha^{-i}\alpha_1\alpha^i,~~~1\leq i,j\leq n-2,~~~
|i-j|\geq 2, \eqno{(25)}
$$
$$
\sigma^{i-1}\sigma_1\sigma^{-(i-1)}\alpha^{-i}\alpha_1\alpha \alpha_1\alpha^{i-1}=
\alpha^{-i}\alpha_1\alpha
\alpha_1\alpha^{i-1}\sigma^{i}\sigma_1\sigma^{-i},~~~2\leq i\leq n-1,
\eqno{(26)}
$$
$$
\sigma^{i-1}\sigma_1\sigma \sigma_1\sigma^{-i}\alpha^{-(i-1)}\alpha_1\alpha^{i-1}=
\alpha^{-i}\alpha_1\alpha^{i}\sigma^{i-1}\sigma_1\sigma \sigma_1\sigma^{-i},~~~2\leq i\leq
n-1.
\eqno{(27)}
$$
New and old generators connecting by formulas
$$
\sigma = \sigma_1\sigma_2\ldots\sigma_{n-1},~~~\alpha = \alpha_{n-1}\alpha_{n-2}\ldots\alpha_1,
~~~ \eqno{(28)}
$$
$$
\sigma_{i+1} = \sigma^i\sigma_1\sigma^{-i},~~~\alpha_{i+1} = \alpha^{-i}\alpha_{1}\alpha^i,
~~~ 1\leq i\leq n-2. \eqno{(29)}
$$
\end{proposition}

\begin{proof} It is well--known [23, Chapter 6] that $B_n$, $n\geq 3$,
is generated by $\sigma_1$ and $\sigma $ and defined by relations (23).
Analogous, symmetric group $S_n$ is generated by $\alpha_1, \alpha $
and defined by relations (24). New and old generators are connecting by formulas (28)
and (29). Hence in presentation of $C_n$ from the Lemma 4.1 we can
replace the generators
$\alpha_i$, $\sigma_i$, $i=1,2,\ldots,n-1$, by generators $\alpha_1$, $\alpha $, $\sigma_1$,
$\sigma $, and replace relations (1)--(2)
by relations (23) and replace relations (13)--(15) by relations (24).
Substitute  in relations (21)--(22) expression $\alpha_i$, $\sigma_j$
through $\alpha $, $\alpha_1$, $\sigma_1$, $\sigma $, we obtain (25)--(27).
This completes the proof.
\end{proof}

Let $V$ be a set of words in an alphabet of variables ranging over
elements of group $G$. The subgroup $V(G)$ of the group $G$
generated by all values of words from $V^{\pm 1}$ is called the {\it
verbal subgroup} defined by the set $V$. The {\it width} wid$(G, V)$ of the
verbal subgroup $V(G)$ with respect to $V$ is defined to be the minimal number $m\in
\mathbb{N} \bigcup \{+\infty \}$ such that any element of $V(G)$
can be presented as a product of $\leq m$ values of words from
$V^{\pm 1}$.

Subgroup is called the {\it proper subgroup} if it is not trivial and
it is not all group. The set $V$ is called the {\it proper set of words},
if the verbal subgroup  $V(F_2)$ is a proper subgroup of $F_2$. If $V$
is not a proper set of words, then it is called
{\it unproper set of words}. The width of verbal subgroup with respect to
unproper set of words is finite [24,
Lemma 1]. If in group  $G$ each finite proper set of words defines the proper verbal
subgroup, then we will said that $G$
{\it is rich by verbal subgroups}. We will said that $G$ {\it has not proper verbal
subgroups of finite width}, if for each finite proper set $V$ the width wid$(G, V)$
is infinite.

Some years ago author has proved  [25] that the braid group $B_n$ for $n\geq 3$
has not proper verbal subgroups of finite width. This result
finishes investigations of some authors  which started with the G.~S.~Makanin's
question from "The Kourovka Notebook" (Unsolved problems
in group theory) [9].

The width of verbal subgroups for HNN--extensions and groups with
one relations were investigated in [26]. The width of derived groups
of some matrix groups over rings were investigated in [27]. The
connection between the width of derived group, commutator length
and equations in free group were studied in [28]. Some Artin groups
were investigated in [29]. In this article were proved

\begin{lemma} ([29, Lemma 3]). If there is an epimorphism from $G$ on
group which has not proper verbal subgroups of finite width and which
is rich by verbal subgroups, then $G$ also has not proper verbal subgroups
of finite width.
\end{lemma}

We will use this lemma and prove

\begin{proposition} The group of bases--conjugating automorphisms
$Cb_n$ for $n\geq 2$ has not proper verbal subgroups of finite width.
\end{proposition}

\begin{proof} By the Theorem 3.1 there exists an epimorphism of $Cb_{n}$ for $n \geq 2$ on $Cb_{2}$.
The group $Cb_2$ is isomorphic to $F_2$. From result of A.~H.~Rhemtulla [24]
we have that $F_2$ has not proper verbal subgroups of finite width and it is evident
that $F_2$ is rich by verbal subgroups.
The result follows from the Lemma 4.2.
\end{proof}

\begin{question}
Is it true that for  $n\geq 2$ the group $C_n$ has not proper verbal subgroups of finite
width?
\end{question}
\vskip 20pt


\section{Expansion of known linear representations of braid group on $C_n$ }

\vskip 12pt

The Burau representation $\psi : B_n\longrightarrow {\rm GL}(W_n)$
maps $B_n$ in group of automorphisms ${\rm GL}(W_n)$ of free $n$ dimension module $W_n$
over ring $\mathbb{Z}[q^{\pm 1}]$ of Laurent polynomials with an integer coefficients.
If $w_1, w_2,\ldots , w_n$ is a free basis of $W_n$,
then $\psi(\sigma_i)$, $i=1,2,\ldots,n-1$ acts by the rule
(we will write $\sigma_i(w_j)$ or $w_j^{\sigma_i}$ instead
$\psi(\sigma_i)(w_j)$ or $w_j^{\psi(\sigma_i)}$ accordingly):
$$
\left.
\begin{array}{ll}
\sigma_i(w_i)=(1-q)w_i+qw_{i+1}, &    \\
\sigma_i(w_{i+1})=w_{i}, &     \\
\sigma_i(w_l)=w_l,~~~ \mbox{if}~~~l\neq i,i+1.&
\end{array}
\right.
$$
We will consider right action of $B_n$ on $W_n$.

A. G. Savushkina [22] proved that group of conjugating
automorphisms $C_n$ generating by automorphisms $\sigma_{1}, \sigma_{2},\ldots, \sigma_{n-1},$
$\alpha_{1}, \alpha_{2}, \ldots, \alpha_{n-1}$ of free group $F_n$.
The set $\sigma_{1}, \sigma_{2},\ldots, \sigma_{n-1}$ generates
braid group $B_n$ and satisfy relations (1)--(2).
The set of automorphisms $\alpha_{1}, \alpha_{2}, \ldots, \alpha_{n-1}$
generates the symmetric group $S_n$ which defines by relations
(13)--(15).
Group $C_n$ is defined by relations of $B_n$, $S_n$ and relations
$$
\alpha_i \sigma_j = \sigma_j \alpha_i~~~\mbox{if}~~~ |i-j|\geq 2,
\eqno{(30)}
$$
$$
\sigma_i\alpha_{i+1}\alpha_i=\alpha_{i+1}\alpha_i\sigma_{i+1}~~~\mbox{if}~~~
i=1,2,...,n-2, \eqno{(31)}
$$
$$
\sigma_{i+1}\sigma_{i}\alpha_{i+1}=\alpha_{i}\sigma_{i+1}\sigma_{i}~~~\mbox{if}~~~
i=1,2,...,n-2. \eqno{(32)}
$$

\begin{proposition} There is a linear representations
$\widetilde{\psi} : C_n\longrightarrow {\rm GL}(W_n)$ of $C_n$ which is an expansion
of Burau representation $\psi : B_n\longrightarrow {\rm GL}(W_n)$.
\end{proposition}

\begin{proof} For the generators of braid group $B_n$ we define
$\widetilde{\psi}(\sigma_i)=\psi(\sigma_i)$. For the generators
$\alpha_{1}, \alpha_{2}, \ldots, \alpha_{n-1}$ define the action on basis
 of $W_n$ by equalities
$$
\left.
\begin{array}{l}
\alpha_i(w_i)=w_{i+1},\\
\alpha_i(w_{i+1})=w_{i},\\
\alpha_i(w_l)=w_l,~~~ \mbox{if}~~~l\neq i,i+1.\\
\end{array}
\right.
$$
Every elements of $C_n$ is a word in alphabet
$\sigma^{\pm 1}_{1}, \sigma^{\pm 1}_{2},\ldots, \sigma^{\pm 1}_{n-1}$,
$\alpha_{1}, \alpha_{2}, \ldots, \alpha_{n-1}$. Hence we defined
the map
$\widetilde{\psi} : C_n\longrightarrow {\rm GL}(W_n)$. We should prove that this
map is homomorphism. For this we should prove that all relations of $C_n$
are carry out in $\widetilde{\psi}(C_n)$.
It is a simple exercise. The proposition is proved.
\end{proof}

We remind (see [14, 15, 16]) definition of faithful linear
representation of Lawrence--Krammer braid group $B_n$. Let $V_n$
be a free module dimension $m=n(n-1)/2$ with basis $v_{ij}$,
$1\leq i < j\leq n$, over ring $\mathbb{Z}[t^{\pm 1}, q^{\pm 1}]$
of Laurent polynomials on two variables. Then representation $\rho
: B_n\longrightarrow {\rm GL}(V_n)$ was defined by action by
$\sigma_i$, $i=1,2,\ldots,n-1$, on module  $V_n$ by equalities
$$
\left.
\begin{array}{l}
\sigma_i(v_{k,i})=(1-q)v_{k,i}+qv_{k,i+1}+q(q-1)v_{i,i+1},\\
 \\
\sigma_i(v_{k,i+1})=v_{k,i}~~~ \mbox{if }~~~k < i,\\
 \\
\sigma_i(v_{i,i+1})=tq^2v_{i,i+1},\\
 \\
\sigma_i(v_{i,l})=tq(q-1)v_{i,i+1}+(1-q)v_{i,l}+qv_{i+1,l}~~~ \mbox{if}~~~i+1 <
l,\\
 \\
\sigma_i(v_{i+1,l})=v_{i,l},\\
 \\
\sigma_i(v_{k,l})=v_{k,l}~~~ \mbox{if }~~~\{ k, l\} \bigcap \{ i,
i+1\}=\emptyset.\\
\end{array}
\right.
$$
From this equalities it is not difficult to find the action by $\sigma_i^{-1}$:
$$
\left.
\begin{array}{l}
\sigma_i^{-1}(v_{k,i})=v_{k,i+1},\\
 \\
\sigma_i^{-1}(v_{k,i+1})=q^{-1}v_{k,i}+q^{-1}(q-1)v_{k,i+1}-t^{-1}q^{-2}(q-1)v_{i,i+1}~~~
 \mbox{if }~~~k < i,\\
  \\
\sigma_i^{-1}(v_{i,i+1})=t^{-1}q^{-2}v_{i,i+1},\\
 \\
\sigma_i^{-1}(v_{i,l})=v_{i+1,l}~~~ \mbox{if }~~~i+1 < l,\\
 \\
\sigma_i^{-1}(v_{i+1,l})=-q^{-2}(q-1)v_{i,i+1}+q^{-1}v_{i,l}+q^{-1}(q-1)v_{i+1,l},\\
 \\
\sigma_i^{-1}(v_{k,l})=v_{k,l}~~~ \mbox{if }~~~\{ k, l\} \bigcap \{ i,
i+1\}=\emptyset.\\
\end{array}
\right.
$$
and the actions by $a_{i,i+1}=\sigma_i^2$ $(i=1,2,\ldots,n-1)$:
$$
\left.
\begin{array}{l}
a_{i,i+1}(v_{k,i})=(q^2-q+1)v_{k,i}+q(1-q)v_{k,i+1}+q(q-1)(tq^2-q+1)v_{i,i+1},\\
 \\
a_{i,i+1}(v_{k,i+1})=(1-q)v_{k,i}+qv_{k,i+1}+q(q-1)v_{i,i+1}~~~ \mbox{if }~~~k <
i,\\
 \\
a_{i,i+1}(v_{i,i+1})=t^2q^4v_{i,i+1},\\
 \\
a_{i,i+1}(v_{i,l})=tq(q-1)(tq^2-q+1)v_{i,i+1}+(q^2-q+1)v_{i,l}+q(1-q)v_{i+1,l}~~~ \mbox{if }~~~i+1 <
l,\\
 \\
a_{i,i+1}(v_{i+1,l})=tq(q-1)v_{i,i+1}+(1-q)v_{i,l}+qv_{i+1,l},\\
 \\
a_{i,i+1}(v_{k,l})=v_{k,l}~~~ \mbox{if }~~~\{ k, l\} \bigcap \{ i,
i+1\}=\emptyset.\\
\end{array}
\right.
$$

We want to expansion the representation of Lawrence--Krammer on
$C_n$. For this we should define the action of $S_n$ on module
$V_n$. Let $S_n$ act on the set $\{1,2,\ldots,n \}$ by
permutation. If $\pi $ is a permutation from $S_n$, then let $\pi
(v_{ij})=v_{\pi(i),\pi(j)}$. Thus automorphism $\alpha_i$ act on
basis of $V_n$ by the rules:
$$
\left.
\begin{array}{l}
\alpha_i(v_{k,i})=v_{k,i+1},\\
 \\
\alpha_i(v_{k,i+1})=v_{k,i}~~~ \mbox{if}~~~k < i,\\
 \\
\alpha_i(v_{i,i+1})=v_{i,i+1},\\
 \\
\alpha_i(v_{i,l})=v_{i+1,l}~~~ \mbox{if}~~~i+1 < l,\\
 \\
\alpha_i(v_{i+1,l})=v_{i,l},\\
 \\
\alpha_i(v_{k,l})=v_{k,l}~~~ \mbox{if}~~~\{ k, l\} \bigcap \{ i,
i+1\}=\emptyset.\\
\end{array}
\right.
$$

It is not difficult to see that this action define a homomorphism from $S_n$
to ${\rm GL}(V_n)$. For this we note that equalities for action of $\alpha_i$ are
obtained from equalities for action of $\sigma_i$ if we take  $q=t=1$.
Hence relations (14)--(15) are valid. The relation $\alpha_i^2=1$ also are valid
that is following from the action of $a_{i,i+1}=\sigma_i^2$ on $V_n$.
Hence we construct the linear representation of $S_n$. As well--known if $n\geq 5$
then $S_n$ has a simple subgroup $A_n$ of index 2 and our
representation are faithful. We have proved

\begin{lemma} The linear representation $S_n\longrightarrow {\rm GL}(V_n)$
constructed above are faithful linear representation of $S_n$.
\end{lemma}

We have the map of $C_n$ in ${\rm GL}(V_n)$ to be defined by
action on generators
 $\alpha_1, \alpha_2,\ldots,\alpha_{n-1}$,
$\sigma_1, \sigma_2,\ldots,\sigma_{n-1}$. This map is keep all relations which include only
generators $\alpha_i$ and keep all relations which include only generators $\sigma_i$.
We should consider remaining relations.
It is not difficult to show that all relations
$$
\alpha_i\sigma_j=\sigma_j\alpha_i,~~ \mbox{if} ~~|i-j|\geq 2,
$$
$$
\sigma_{i+1}\sigma_{i}\alpha_{i+1}=\alpha_i\sigma_{i+1}\sigma_i,
$$
are keep.

Consider relation $\sigma_i\alpha_{i+1}\alpha_i=\alpha_{i+1}\alpha_i\sigma_{i+1}$.
Acting by the left--hand side of this relations on basis of $V_n$
we have
$$
\left.
\begin{array}{l}
v_{k,i}^{\sigma_i\alpha_{i+1}\alpha_i}=(1-q)v_{k,i+1}+qv_{k,i+2}+q(q-1)v_{i+1,i+2},\\
 \\
v_{k,i+1}^{\sigma_i\alpha_{i+1}\alpha_i}=v_{k,i+1}~~~ \mbox{if}~~~k <
i,\\
 \\
v_{i,i+1}^{\sigma_i\alpha_{i+1}\alpha_i}=tq^2v_{i+1,i+2},\\
 \\
v_{i,l}^{\sigma_i\alpha_{i+1}\alpha_i}=tq(q-1)v_{i+1,i+2}+(1-q)v_{i+1,l}+qv_{i+2,l}~~~ \mbox{if}~~~i+2 <
l,\\
 \\
v_{i,i+2}^{\sigma_i\alpha_{i+1}\alpha_i}=tq(q-1)v_{i+1,i+2}+(1-q)v_{i,i+1}+qv_{i,i+2},\\
 \\
v_{i+1,l}^{\sigma_i\alpha_{i+1}\alpha_i}=v_{i+1,l}~~~ \mbox{if}~~~l >
i+2,\\
 \\
v_{i+1,i+2}^{\sigma_i\alpha_{i+1}\alpha_i}=v_{i,i+1},\\
 \\
v_{k,l}^{\sigma_i\alpha_{i+1}\alpha_i}=v_{k,l}~~~ \mbox{if}~~~\{ k, l\} \bigcap \{ i,
i+1, i+2\}=\emptyset,\\
 \\
v_{k,i+2}^{\sigma_i\alpha_{i+1}\alpha_i}=v_{k,i}~~~ \mbox{if}~~~k <
i,\\
 \\
v_{i+2,l}^{\sigma_i\alpha_{i+1}\alpha_i}=v_{i,l}.\\
\end{array}
\right.
$$
Acting by the right--hand side of this relations, we have
$$
\left.
\begin{array}{l}
v_{k,i}^{\alpha_{i+1}\alpha_i\sigma_{i+1}}=(1-q)v_{k,i+1}+qv_{k,i+2}+q(q-1)v_{i+1,i+2},\\
 \\
v_{k,i+1}^{\alpha_{i+1}\alpha_i\sigma_{i+1}}=v_{k,i+1}~~~ \mbox{if}~~~k <
i,\\
 \\
v_{i,i+1}^{\alpha_{i+1}\alpha_i\sigma_{i+1}}=tq^2v_{i+1,i+2},\\
 \\
v_{i,l}^{\alpha_{i+1}\alpha_i\sigma_{i+1}}=tq(q-1)v_{i+1,i+2}+(1-q)v_{i+1,l}+qv_{i+2,l}~~~ \mbox{if}~~~i+2 <
l,\\
 \\
v_{i,i+2}^{\alpha_{i+1}\alpha_i\sigma_{i+1}}=(1-q)v_{i,i+1}+qv_{i,i+2}+q(q-1)v_{i+1,i+2},\\
 \\
v_{i+1,l}^{\alpha_{i+1}\alpha_i\sigma_{i+1}}=v_{i+1,l}~~~ \mbox{if}~~~l >
i+2,\\
 \\
v_{i+1,i+2}^{\alpha_{i+1}\alpha_i\sigma_{i+1}}=v_{i,i+1},\\
 \\
v_{k,l}^{\alpha_{i+1}\alpha_i\sigma_{i+1}}=v_{k,l}~~~ \mbox{if}~~~\{ k, l\} \bigcap \{ i,
i+1, i+2\}=\emptyset,\\
 \\
v_{k,i+2}^{\alpha_{i+1}\alpha_i\sigma_{i+1}}=v_{k,i}~~~ \mbox{if}~~~k <
i,\\
 \\
v_{i+2,l}^{\alpha_{i+1}\alpha_i\sigma_{i+1}}=v_{i,l}.\\
\end{array}
\right.
$$

Comparing the actions of automorphisms $\sigma_i\alpha_{i+1}\alpha_i$
and
$\alpha_{i+1}\alpha_i\sigma_{i+1}$ on $v_{i,i+2}$ we see that they are equal if and only if $t=1$.
Therefore the following proposition is true

\begin{proposition} The map
$\rho : C_n\longrightarrow {\rm GL}(V_n)$
constructed above is a linear representation if and only if $t=1$.
\end{proposition}

For $t=1$ $B_n$--module $V_n$ is isomorphic to the symmetric square of $B_n$--module $W$ and
in this case representation $\rho : C_n\longrightarrow {\rm GL}(V_n)$
is unfaithful for
$n\geq 5$. The question of the faithfulity of this representation for $n=3,4$
remained open.

We can define the action of $\alpha_i$
on $V_n$ by another way. Not difficult to prove

\begin{lemma} For every nonzero number  $x\in \mathbb{C}$ the map
$\rho_1 : S_n\longrightarrow {\rm GL}(V_n)$ defining by action of
 $\alpha_i$ on the $V_n$ by equalities
$$
\left.
\begin{array}{l}
\alpha_i(v_{k,i})=xv_{k,i+1},\\
 \\
\alpha_i(v_{k,i+1})=x^{-1}v_{k,i}~~~ \mbox{if}~~~k < i,\\
 \\
\alpha_i(v_{i,i+1})=v_{i,i+1},\\
 \\
\alpha_i(v_{i,l})=xv_{i+1,l}~~~ \mbox{if}~~~i+1 < l,\\
 \\
\alpha_i(v_{i+1,l})=x^{-1}v_{i,l},\\
 \\
\alpha_i(v_{k,l})=v_{k,l}~~~ \mbox{if}~~~\{ k, l\} \bigcap \{ i,
i+1\}=\emptyset.\\
\end{array}
\right.
$$
is the faithful linear representation.
\end{lemma}

If we will expand the map $\rho_1$ on all $C_n$ defining $\rho_1(\sigma_i)=\rho(\sigma_i)$,
then we will have the map $\rho_1 : C_n\longrightarrow {\rm GL}(V_n)$.
We see that this map keeps relations
$$
\sigma_i\alpha_{i+1}\alpha_i=\alpha_{i+1}\alpha_i\sigma_{i+1},
$$
$$
\sigma_{i+1}\sigma_{i}\alpha_{i+1}=\alpha_i\sigma_{i+1}\sigma_i
$$
if we take $x=t^{-1/3}$. but in this case the relation
$\alpha_i\sigma_j=\sigma_j\alpha_i$ is not keep since
$$
v_{i,j}^{\alpha_i\sigma_j}=x(1-q)v_{i+1,j}+xqv_{i+1,j+1}+xq(q-1)v_{j.j+1},
$$
$$
v_{i,j}^{\sigma_j\alpha_i}=x(1-q)v_{i+1,j}+xqv_{i+1,j+1}+q(q-1)v_{j.j+1},
$$
and $v_{i,j}^{\alpha_i\sigma_j}=v_{i,j}^{\sigma_j\alpha_i}$
if and only if $x=1$. Since there  not exist relations
$\alpha_i\sigma_j=\sigma_j\alpha_i$ in
$C_3$ then we have

\begin{proposition} The linear representation $\rho_1 : C_3\longrightarrow {\rm GL}(V_3)$
is expansion of Lawrence--Krammer representation on $C_3$.
\end{proposition}

\vskip 20pt

\begin{center}

\section{Linear representations of braid groups of some manifolds}

\end{center}

\vskip 5pt

As it was proved  (see [7, Theorem 1.11; 32])
that braid group $B_n(S^2)$ of the 2--sphere $S^2$ admits a
presentation with generators
$\delta_1, \delta_2,\ldots,\delta_{n-1}$
and defining relations:
$$
\delta_i\delta_{i+1}\delta_{i}=\delta_{i+1}\delta_{i}\delta_{i+1}~~~\mbox{if}~~~
i=1,2,...,n-2,
$$
$$
\delta_i\delta_{j}=\delta_{j}\delta_{i}~~~\mbox{if}~~~ |i-j|~~~ \geq
2,
$$
$$
\delta_1\delta_{2}\ldots \delta_{n-2}\delta_{n-1}^2\delta_{n-2}\ldots
\delta_{2}\delta_1=1.
$$
From this relations we see that $B_n(S^2)$ is a homomorphic image
of
$B_n$. Since $B_2(S^2)$ is a cyclic group of order 2 and
$B_3(S^2)$ is a metacyclic group of order 12, then we will consider $n>3$.

R.~Gillette and J.~Van Buskirk [30] have studied the structure of $B_n(S^2)$.
We recall some of their results. Let
$$
a_{i,i}=1,~~~a_{i,j}=\delta_i^{-1}\delta_{i+1}^{-1}\ldots \delta_{j-2}^{-1}\delta_{j-1}^2\delta_{j-2}\ldots
\delta_{i+1}\delta_i,~~~1\leq i< j\leq n. \eqno{(33)}
$$
These elements generate the pure braid group $P_n(S^2)$ which is a kernel
of homomorphism $\nu : B_n(S^2)\longrightarrow S_n$ and $\nu $
send
$\delta_i$ in transposition $(i,i+1)$,
$1\leq i\leq n-1$, from $S_n$.  Let define the subgroup
$A_{n-i+1}=<a_{i,i+1}, a_{i,i+2},\ldots,a_{i,n}>$ for each
$i=1,2,\ldots,n-1$.
Subgroup $A_n$
is normal in $P_n(S^2)$ and we have the short exact sequence
$$
1\longrightarrow A_n\longrightarrow P_n(S^2)\longrightarrow
P_{n-1}(S^2)\longrightarrow 1
$$
and $P_n(S^2)$ is a semi direct product: $P_n(S^2)=A_n\leftthreetimes
P_{n-1}(S^2)$.
The generators of $A_n$ are connected by the relation
$$
a_{1,2}a_{1,3}\ldots a_{1,n}=1.
$$
Since the following relations are hold in $P_n(S^2)$
$$
a_{i,i+1}a_{i,i+2}\ldots a_{i,i+n-1}=1,~~~i=1,2,\ldots,n-1,
$$
where $a_{j,i}=a_{i,j}$ if $j > i$ and indices are taken mod$n$,
then
$$
a_{i,n}=a_{i,n-1}^{-1}\ldots a_{i,i+1}^{-1}a_{i-1,i}^{-1}\ldots
a_{1,i}^{-1},~~~i=1,2,\ldots,n-1.
$$
Using this formulas we can exclude  $a_{1,n}, a_{2,n},\ldots, a_{n-1,n}$ from the set
of generators of $P_n(S^2)$.
Hence $A_{n-i+1}$ is freely generated by
$a_{i,i+1}, a_{i,i+2},\ldots, a_{i,n-1}$.

The center of $B_n(S^2)$ is the cyclic group of order 2 generated
by the Dirac braid
$$
\Delta_n=(\delta_1\delta_2 \ldots \delta_{n-2})^{n-1}=
(a_{1,2}a_{1,3}\ldots a_{1,n-1})(a_{2,3}a_{2,4}\ldots a_{2,n-1})\ldots
(a_{n-2,n-1}).
$$

The group $P_n(S^2)$ was split in the semi direct product
$$
P_n(S^2)=A_n\leftthreetimes (A_{n-1}\leftthreetimes (\ldots
\leftthreetimes A_3)\ldots).
$$
Let $L_n=A_n\leftthreetimes (A_{n-1}\leftthreetimes (\ldots
\leftthreetimes A_4)\ldots)$. Since $A_3$ is a cyclic group
generating by
$a_{n-2,n-1}$, then $P_n(S^2)=L_n\times <\Delta_n>\simeq L_n\times \mathbb{Z}_2$.
Group $L_n$ is isomorphic to subgroup
$U_n\leftthreetimes (U_{n-1}\leftthreetimes (\ldots
\leftthreetimes U_4)\ldots)\leq P_n$.

 The braid group of sphere closed relates with the group $M(0,n)$
of mapping classes of the $n$--punctured 2--sphere. Presentation of the $M(0,n)$
derived from that of $B_n(S^2)$ by adding the relation:
$\Delta_n=(\delta_1\delta_2 \ldots \delta_{n-2})^{n-1}=1$. There is an epimorphism from
 $B_n(S^2)$ to $M(0,n)$. The kernel of this epimorphism is equal
 to the center of
$B_n(S^2)$.
A mapping class group of the $n$--punctured 2--sphere can be visualized
as an $n$--string braid between concentric 2--spheres, where the
inner sphere is free to execute full revolutions.

For $M(0,n)$ there is a short exact sequence
$$
1\longrightarrow L_n\longrightarrow M(0,n)\longrightarrow
S_{n}\longrightarrow 1,
$$
where $L_n$ is a kernel of epimorphism $\nu : M(0,n)\longrightarrow S_n$ that map
 $\delta_i$ to the transposition $(i,i+1)$, $i=1,2,\ldots, n-2$.

A. I. Malcev [31, Lemma 1; 21] have proved that if $H$ is subgroup
of finite index in group $G$ and $H$ is linear, then $G$ is linear
too. Recall his construction. Let $\vert G : H \vert =m$ and $\psi
: H\longrightarrow {\rm GL}_{l}(F)$ is a faithful linear
representation of $H$ by matrices of order $l$ over field $F$. In
this case $G$ is a union of right cosets
$$
G=He \bigcup Hg_2 \bigcup \ldots \bigcup Hg_m.
$$
For each $g\in G$ we can write the product $g_ig$ uniquely in the form
 $h_ig_{n_i}$, $h_i \in H$. Therefore for each $g$ we have the sequence $h_1, h_2, \ldots, h_m$
from $H$ and the sequence numbers $n_1, n_2, \ldots, n_m$. For the
sequence $\{ n_i \}$ we construct the matrix $D(n_i)$ with integer
coefficients by the rule:
$$
D(n_i)=\parallel d_{j,k} \parallel \in
M_{lm}(\mathbb{Z}),~~~d_{j,n_j}=E_l,~~~d_{j,k}=0~~~
\mbox{if}~~~k\neq n_j,
$$
where $E_l$ is identical matrix of order $l$. Then we can map $g\in G$
to the matrix
$$
{\rm diag}(\psi(h_1), \psi(h_2),\ldots, \psi(h_m))D(n_i).
$$
Hence we constructed the linear representation of $G$ in ${\rm GL}_{lm}(F)$.
This linear representation is faithful. Using this construction we
will prove

\begin{theorem} Let $R=\mathbb{Z}[t^{\pm 1}, q^{\pm 1}]$ be a ring
of Laurent polynomials on two variables.
Then the groups $M(0,n)$ and $B_n(S^2)$ are linear for each $n\geq 2$.
For
each $n\geq 4$ there are inclusion maps $\varphi : M(0,n)\longrightarrow {\rm GL}_m(R)$
and
$\varphi_1 : B_n(S^2)\longrightarrow {\rm GL}_{m_1}(R)$, where $m=(n-1)(n-2)n!/2$, $m_1=2m$.
\end{theorem}

\begin{proof}
As we noted above, if $n=2,3$ groups $M(0,n)$, $B_n(S^2)$
are finite and hence are linear. Since $L_n$ is isomorphic to subgroup of $B_{n}$,
then there are faithful linear representation $\rho : L_n\longrightarrow {\rm GL}_l(R)$ for
 $l=(n-1)(n-2)/2$. This representation is induced by Lawrence--Krammer representation of $B_n$.

Let $m_1, m_2, \ldots, m_{n!}$ be coset representatives of $M(0,n)$ by
subgroup $L_n$. Since $M(0,n)$ is generated by
$\delta_1, \delta_2, \ldots, \delta_{n-1}$ then we can define
$\varphi $ only on these elements. Each generator $\delta_k$
acts on set of coset representatives by permutation. We find
$$
m_i\delta_k =
h_i^km_{\pi_k(i)},~~~i=1,2,\ldots,n!,~~~k=1,2,\ldots,n-1,
$$
where symbol $k$ in the upper part is an index but not exponent (this rule will be true to the end
of this proof). To compare with  $\delta_k$ the matrix
$$
\varphi(\delta_k) = {\rm diag}(\rho (h_1^k), \rho (h_2^k),\ldots, \rho (h_{n!}^k))
\pi(\delta_k)\in {\rm GL}_m(R),
$$
where ${\rm diag}(\rho (h_1^k), \rho (h_2^k),\ldots, \rho (h_{n!}^k))$
is a block--diagonal matrix; $\pi(\delta_k)$ is a block--monomial matrix in which
the block on position $(j,\pi(j))$ is an identical matrix of order $l$, but the block on position
 $(j,s)$ for $s\neq \pi(j)$
is a null matrix of order $l$. This linear representation $\varphi $
is a faithful representation of
$M(0,n)$.

Consider the group $B_n(S^2)$. This group contains the linear subgroup $L_n$
index $2n!$. As a set of right coset representatives of $B_n(S^2)$ by subgroup $L_n$
we take elements $m_i\Delta_n^{\epsilon }$,
$i=1,2,\ldots,n!$, $\epsilon =0,1,$ where $\Delta_n$ is the generator of center of $B_n(S^2)$.
Since $B_n(S^2)$ is generated by $\delta_1, \delta_2, \ldots, \delta_{n-1}$ then
it is enough to define representation $\varphi_1$ on these elements. Let us ordere the coset representatives
 $m_i\Delta_n^{\epsilon }$ and denote them by
 $n_1, n_2,\ldots, n_{2n!}$. Each generator $\delta_k$ of
$B_n(S^2)$
acts on these representatives by the rules
$$
n_j\delta_k =
g_j^k n_{\pi_k(j)},~~~j=1,2,\ldots,2n!,~~~k=1,2,\ldots,n-1.
$$
To compare with $\delta_k$ the matrix
$$
\varphi_1(\delta_k) = {\rm diag}(\rho (g_1^k), \rho (g_2^k),\ldots, \rho (g_{2n!}^k))
\pi(\delta_k)\in {\rm GL}_{2m}(R),
$$
where ${\rm diag}(\rho (g_1^k), \rho (g_2^k),\ldots, \rho (g_{2n!}^k))$
is a block--diagonal matrix and $\pi(\delta_k)$ is a block--monomial matrix in which the block
on position $(j,\pi(j))$ is an identical matrix of order $l$ but the block on position $(j,s)$ for
 $s\neq \pi(j)$
is a null matrix of order $l$. This linear representation $\varphi_1$
is a faithful representation of $B_n(S^2)$ over ring $R$. Since $R$
is included in the field of complex numbers $\mathbb{C}$ (it is enough to take numbers $t$ and $q$
which are non--zero transcendental over $\mathbb{Q}$) then we obtain the required assertions.
\end{proof}

\vskip 5pt

Let $P^2$ be a projective plane. The braid group $B_n(P^2)$,
$n\geq 1,$ admits a presentation [32] with generators
$\delta_1, \delta_2,\ldots,\delta_{n-1},$ $\rho_1, \rho_2,\ldots,\rho_{n}$ and
defining relations:
$$
\delta_i\delta_{i+1}\delta_{i}=\delta_{i+1}\delta_{i}\delta_{i+1}~~~\mbox{if}~~~ i=1,2,...,n-2,
$$
$$
\delta_i\delta_{j}=\delta_{j}\delta_{i}~~~\mbox{if}~~~ |i-j| \geq 2,
$$
$$
\delta_i\rho_{j}=\rho_{j}\delta_{i},~~~j\neq i,i+1,
$$
$$
\rho_i = \delta_i\rho_{i+1}\delta_{i}~~~\mbox{if}~~~ i=1,2,...,n-1,
$$
$$
\rho_{i+1}^{-1}\rho_{i}^{-1}\rho_{i+1}\rho_{i}=\delta_{i}^2~~~\mbox{if}~~~ i=1,2,...,n-1,
$$
$$
\delta_1\delta_{2}\ldots \delta_{n-2}\delta_{n-1}^2\delta_{n-2}\ldots
\delta_{2}\delta_1=\rho_1^2.
$$

It is clear that the map taking  $\delta_i$ onto the transposition
$(i,i+1)$, $i=1, 2, \ldots, n-1,$
and each $\rho_j$ onto the identity is  homomorphism of $B_n(P^2)$
onto the symmetric group $S_n$, with kernel $P_n(P^2)$.

The group $B_1(P^2)$ is a cyclic group of order 2,  $B_2(P^2)$ is a finite group of order 16.
For $n\geq 3$ the group $B_n(P^2)$ is infinite. Consider the case $n=3$. The group $P_3(P^2)$
contains a free subgroup $A_3(P^2)=<\rho_1, a_2, a_3 \parallel a_2a_3=\rho_1^2>$,
where $a_2=\delta_1^2$, $a_3=\delta_1^{-1}\delta_2^2\delta_1$. The
factor group $P_3(P^2)/A_3(P^2)\simeq P_2(P^2)$ is a quaternion group of order 8.
The set
$$
\{ 1, \delta_1,
\delta_2, \delta_2\delta_1, \delta_1\delta_2, \delta_1\delta_2\delta_1 \}
$$
is a right coset representatives of $B_3(P^2)$ by $P_3(P^2)$.
Note that $P_2(P^2)$ has a presentation
$$
P_2(P^2) = <\rho_2, \rho_3 \parallel
\rho_2^2=\rho_3^2=(\rho_3\rho_2)^2>.
$$
In $P_3(P^2)$ the following relations are true:
$$
\rho_3^{-\epsilon }a_2\rho_3^{\epsilon
}=a_2,~~~\epsilon=\pm 1,~~~\rho_j^{-1}\rho_1\rho_j=\rho_1a_j^{-1},~~j=2,3,
$$
$$
\rho_2^{-1}a_3\rho_2=a_2a_3a_2^{-1},~~~\rho_2^{-1}a_2\rho_2=\rho_1a_2^{-1}\rho_1^{-1}.
$$
From defining relations of $B_3(P^2)$ we get $\delta_1$, $\delta_2$.

\begin{lemma} In $B_3(P^2)$ the following rules of conjugating
are true
$$
a_2^{\delta_1}=a_2,~~~\rho_1^{\delta_1}=\rho_2a_2,~~~\rho_2^{\delta_1}=a_2^{-1}\rho_1,
~~~\rho_3^{\delta_1}=\rho_3,
$$
$$
a_2^{\delta_2}=\rho_1^2a_2^{-1},~~~\rho_1^{\delta_2}=\rho_1,
~~~\rho_2^{\delta_2}=\rho_3\delta_2^2,~~~\rho_3^{\delta_2}=\delta_2^{-2}\rho_2.
$$
\end{lemma}

Since $A_3(P^2)$ is a free group with free generators $\rho_1$, $a_2$
 then using Sanov's representation
[21] we can include $A_3(P^2)$ in ${\rm SL}_2(\mathbb{Z})$
assuming
$$
\overline{\rho}_1=
\left(
\begin{array}{cc}
1 &  2  \\
0 &  1   \\
\end{array}
\right),~~~
\overline{a}_2=
\left(
\begin{array}{cc}
1 &  0  \\
2 &  1   \\
\end{array}
\right)~~~
$$
as images of $\rho_1$ and $a_2$ accordingly.

The group $P_3(P^2)$ is an extension of $A_3(P^2)$ by the group of quaternion of order 8.
Let as chose the set
$$
\{ e, \rho_2, \rho^2_2, \rho^3_2, \rho_3, \rho_3\rho_2, \rho_3\rho^2_2, \rho_3\rho^3_2
\}
$$
as the set of coset representatives of $P_3(P^2)$ by $A_3(P^2)$.
Thus the free group $A_3(P^2)$ has index 48 in
$B_3(P^2)$. Using Malcev's construction we get

\begin{theorem} Braid group $B_3(P^2)$ of projective plane $P^2$ is included in
 ${\rm SL}_{96}(\mathbb{Z})$.
\end{theorem}

\vskip 20pt

\begin{center}

\section{The faithful linear representation of ${\rm Aut}(F_2)$}

\end{center}

\vskip 5pt

It is known [17] that ${\rm Aut}(F_n)$ is not linear for $n\geq 3$. If $n=2$ then in [20]
was proved that ${\rm Aut}(F_2)$ is linear if and only if
$B_4$ is. Since $B_4$ is linear then we can construct the faithful linear representation
of ${\rm Aut}(F_2)$.
Recall some results from [20].

Let $B_4^{*}=B_4/Z(B_4)$ where $Z(B_4)$ is a center of $B_4$. As
it
was noted above
$Z(B_4)$ is an infinite cyclic group generating by
$\Delta_4=(\sigma_1\sigma_2\sigma_3)^4$.
The group $B_4^{*}$ is isomorphic to ${\rm Aut}^+(F_2)$ which is
preimage of ${\rm SL}_2(\mathbb{Z})$ by homomorphism
$\xi : {\rm Aut}(F_2)\longrightarrow {\rm GL}_2(\mathbb{Z})$ sending automorphism of $F_2$
to automorphism of $\mathbb{Z}^2$. The group ${\rm
Aut}^+(F_2)$ is subgroup of index 2 in ${\rm Aut}(F_2)$.

Let $\rho : B_4\longrightarrow {\rm GL}_6(\mathbb{C})$ be a faithful linear representation of
Lawrense--Krammer. M.~Zino [33] have proved that this representation is irreducible.

Denote by $F$ subgroup of $B_4$ generating by
$x=\sigma_1\sigma^{-1}_3$, $y=\sigma_2\sigma_1\sigma^{-1}_3\sigma^{-1}_2$.
This group is free with free generators
 $x$, $y$ and is normal in $B_4$. Inner automorphisms of $B_4$
 induced automorphisms of $F$, i.~e.,
there is epimorphism $h : B_4\longrightarrow {\rm Aut}^+(F_2)$.
The kernel of this epimorphism is equal to the center of
$B_4$ and images of generators $\sigma_1$, $\sigma_2$, $\sigma_3$
defining automorphisms
$$
h(\sigma_{1})=\alpha_1 : \left\{
\begin{array}{ll}
x \longmapsto x, &  \\
y \longmapsto yx^{-1}, & \\
\end{array} \right.~~~
h(\sigma_{2})=\alpha_2 : \left\{
\begin{array}{ll}
x \longmapsto y, &  \\
y \longmapsto yx^{-1}y, & \\
\end{array} \right.~~~
h(\sigma_{3})=\alpha_3 : \left\{
\begin{array}{ll}
x \longmapsto x, &  \\
y \longmapsto x^{-1}y. & \\
\end{array} \right.~~~
$$

The whole group ${\rm Aut}(F_2)$ is generated by automorphisms
$$
P : \left\{
\begin{array}{ll}
x \longmapsto y, &  \\
y \longmapsto x, & \\
\end{array} \right.~~~
\omega : \left\{
\begin{array}{ll}
x \longmapsto x^{-1}, &  \\
y \longmapsto y, & \\
\end{array} \right.~~~
U : \left\{
\begin{array}{ll}
x \longmapsto xy, &  \\
y \longmapsto y, & \\
\end{array} \right.~~~
$$
and defined by relations
$$
P^2=\omega^2=(\omega P)^4=(P\omega PU)^2=(UP\omega
)^3=[\omega, \omega U\omega ]=1.
$$
The following formulas are hold
$$
\alpha_1=PU^{-1}P,~~~\alpha_2=PU\omega U^{-1},~~~\alpha_3=P\omega U\omega P.
$$

Using the Lawrense--Krammer representation $\rho : B_4\longrightarrow {\rm GL}_6(\mathbb{C})$
we find the images of generators of $B_4$:
$$
\rho_1=\rho(\sigma_1)=
\left(
\begin{array}{cccccc}
tq^2 &  0 & 0 & 0 & 0 & 0  \\
tq(q-1) &  1-q & 0 & q & 0 & 0   \\
tq(q-1) & 0 & 1-q & 0 & q & 0    \\
0 &  1 & 0 & 0 & 0 & 0   \\
0 &  0 & 1 & 0 & 0 & 0   \\
0 &  0 & 0 & 0 & 0 & 1   \\
\end{array}
\right),
$$
$$
\rho_2=\rho(\sigma_2)=
\left(
\begin{array}{cccccc}
1-q &  q & 0 & q(q-1) & 0 & 0  \\
1 &  0 & 0 & 0 & 0 & 0   \\
0 & 0 & 1 & 0 & 0 & 0    \\
0 &  0 & 0 & tq^2 & 0 & 0   \\
0 &  0 & 0 & tq(q-1) & 1-q & q   \\
0 &  0 & 0 & 0 & 1 & 0   \\
\end{array}
\right),
$$
$$
\rho_3=\rho(\sigma_3)=
\left(
\begin{array}{cccccc}
1 &  0 & 0 & 0 & 0 & 0  \\
0 &  1-q & q & 0 & 0 & q(q-1)   \\
0 & 1 & 0 & 0 & 0 & 0    \\
0 &  0 & 0 & 1-q & q & q(q-1)   \\
0 &  0 & 1 & 0 & 0 & 0   \\
0 &  0 & 0 & 0 & 0 & tq^2   \\
\end{array}
\right).
$$
It is not difficult to check that $(\rho_1\rho_2\rho_3)^4=t^2q^8E_6$ is the image of center of $B_4$.
Let $\mu=1/\sqrt[12]{t^2q^8}$ and define $\overline{\rho}(\sigma_i)=\mu \rho(\sigma_i)=\mu \rho_i$.
Then $\overline{\rho}$ induce the faithful linear representation of $B_4^{*}\simeq {\rm Aut}^+F_2$.

Since ${\rm Aut}^+F_2$ is subgroup of index 2 in ${\rm Aut}F_2$
then  $e$ and
$\omega $ is a right coset representatieves. It is not difficult to
check

\begin{lemma} In ${\rm Aut}F_2$  the following equalities are hold
$$
\begin{array}{ll}
1) & \omega \alpha_1=\alpha^{-1}_1\omega,\\
2) & \omega
\alpha_2=\alpha_3\alpha_1\alpha^{-1}_2\alpha^{-1}_1\alpha^{-1}_3\omega,\\
3) & \omega \alpha_3=\alpha^{-1}_3\omega.
\end{array}
$$
\end{lemma}

Since ${\rm Aut}F_2$ is generated by automorphisms $\alpha_1, $ $\alpha_2,$ $\alpha_3$,
$\omega,$ then acting on coset representatives of ${\rm Aut}F_2$ by ${\rm Aut}^+F_2$
we find the correspondence diagonal and monomial matrix. Define
$$
\psi(\alpha_1)={\rm diag}(\overline{\rho }_1, (\overline{\rho
}_1)^{-1}),~~~\psi(\alpha_2)={\rm diag}(\overline{\rho }_2,
\overline{\rho }_3\overline{\rho }_1(\overline{\rho }_2)^{-1}(\overline{\rho
}_1)^{-1}(\overline{\rho}_3)^{-1}),
$$
$$
\psi(\alpha_3)={\rm diag}(\overline{\rho }_3, (\overline{\rho
}_3)^{-1}),~~~
\psi(\omega )=
\left(
\begin{array}{cc}
0 &  E_6  \\
E_6 &  0    \\
\end{array}
\right),
$$
where $E_6$ is the identical matrix of order 6 and $\overline{\rho }_i$
defined above. If we conjugate this representation by matrix
$c={\rm diag}(E_6, \overline{\rho }_3\overline{\rho}_1)$ we receive the new representation:
$$
\overline{\psi }(\alpha_i)=c^{-1}\psi(\alpha_i)c=
\left(
\begin{array}{cc}
\overline{\rho }_i &  0  \\
0 &  \overline{\rho }_i^{-1}    \\
\end{array}
\right),~~~i=1,2,3;~~~
\overline{\psi }(\omega )=c^{-1}\psi(\omega )c=
\left(
\begin{array}{cc}
0 & \overline{\rho }_3 \overline{\rho }_1  \\
(\overline{\rho }_3 \overline{\rho }_1)^{-1} &  0    \\
\end{array}
\right).
\eqno{(34)}
$$

The following theorem is true

\begin{theorem} The map
$\overline{\psi } : {\rm Aut}F_2\longrightarrow {\rm GL}_{12}(\mathbb{Z}[t^{\pm 1}, q^{\pm 1}])$,
defining on generators in (34)
is a faithful linear representation of  ${\rm Aut}F_2$.
\end{theorem}

\vskip 30pt

\centerline{REFERENCES} \vskip 12pt
\begin{enumerate}
\item
    D.~Z.~Djokovic, The structure of the automorphism group of a free
group on two generators, Proc. Amer. Math. Soc., 88, N~2 (1983),
218--220.
\item
    G.~T.~Kozlov, The structure of the automorphism group ${\rm Aut}(F_2)$, Algebra, Logika i
prilozhenija, Irkutsk, IGU, 1994, 28--32 (in Russian).
\item
    R.~C.~Lyndon and P.~E.~Schupp, Combinatorial group theory, Springer--Verlag, 1977.
\item
    S.~Krstic, J.~McCool, The non-finite presentability of ${\rm IA}(F_3)$ and
${\rm GL}_2(\mathbb{Z}[t,t^{-1}])$, Inven. Math., 129, N~3
(1997), 595--606.
\item
    J.~McCool, On basis--conjugating automorphisms of free groups, Can. J. Math., 38, N~6 (1986), 1525--1529.
\item
    A.~A.~Markov, Foundations of the algebraic theory of braids,
Trudy Math. Inst. Steklov, 1945, 16, 1--54 (in Russian).
\item
    J. S. Birman, Braids, links and mapping class group,
Princeton--Tokyo: Univ. press, 1974.
\item
    D.~J.~Collins, N.~D.~Gilbert, Structure and torsion in automorphism groups of free products,
Quart. J. Math. Oxford, 41, N~162 (1990), 155--178.
\item
   The Kourovka Notebook (Unsolved problems
in group theory), 15th edn., Institute of Mathematics SO RAN, Novosibirsk 2002.
\item
    W. Burau, Uber Zopfgruppen und gleichsinnig verdrillte
Verkettungen, Abh. Math. Semin. Hamburg Univ, 11, 1936,
179--186.
\item
    J. A. Moody, The Burau representation of the braid group
    $B_n$ is unfaithful for large $n$, Bull. Amer. Math. Soc.
25, N~2 (1991), 379-384.
\item
    D. D. Long, M. Paton, The Burau representation
is not faithful for $n\geq 6$, Topology, 32, N~2 (1993), 439-447.
\item
    S. Bigelow, The Burau representation of $B_5$
is not faithful, Geom. Topology, 3,  (1999), 397-404.
\item
    R. J. Lawrence, Homological representation
of the Hecke Algebra, Commun. Math. Phys., 135, N~1 (1990), 141-191.
\item
    D.~Krammer, Braid groups are linear, Annals of Math.,
155, N~1 (2002), 131-156.
\item
    S.~Bigelow, Braid groups are linear, J. Amer.  Math. Soc.,
14, (2001), 471-486.
\item
E. Formanek, C. Procesi, The automorphism groups of a free group
is not linear, J. Algebra, 149, N~2 (1992), 494--499.
\item
    V.~G.~Bardakov, M.~V.~Nechshadim, Some property of braid
groups of compact oriented 2--manifolds, IV International
Algebraic Conference dedicated to the 60th anniversary of
Yu.~I.~Merzljakov, Novosibirsk, Institute of Mathematics, 2000,
9--13 (in Russian).
\item
    S.~Bigelow, R.~D.~Budney, The mapping class group of genus two
    surface is linear, to appear in "Algebraic and Geometric
    Topology".
\item
J. L. Dyer, E. Formanek, E. K. Grossman, On the linearity of
automorphism groups of free groups, Arch. Math., 38, N~5 (1982),
404--409.
\item
M.~I.~Kargapolov, Yu.~I.~Merzljakov, Fundamentals of the theory of
groups, New York: Springer, 1979.
\item
    A.~G.~Savushkina, On group of conjugating automorphisms of free group,
Matem. Zametki,
60, N~1 (1996), 92--108 (in Russian).
\item
H.~S.~M.~Coxeter, H.~W.~O.~J.~Mozer, Generators and relations for
discrete
groups, Berlin; Heidelberg; New York: Springer, 1957.
\item
        A.~H.~Rhemtulla, A problem of bounded expressibility in free
products, Proc. Cambridge Phil. Soc., 64, N~3 (1969), 573--584.
\item
V. G. Bardakov, On the theory of braid groups, Russ. Acad. Sci., Sb. Math., 76, N~1 (1993),
123--153 (Zbl 0798.20029).
\item
V. G. Bardakov, On the width of verbal subgroups of certain free constructions,
Algebra Logic, 36, N~5 (1997),
288--301 (Zbl 0941.20017).
\item
V. G. Bardakov, On the decomposition of automorphisms of free modules into simple factors,
Izv. Ross. Acad. Nauk, Ser. Mat., 59, N~2 (1995),
109--128 (Zbl 0896.20031).
\item
V. G. Bardakov, Computing the commutator length in free groups,
Algebra Log., 39, N~4 (2000), 395--440 (Zbl 096.20019).
\item
V. G. Bardakov, The width of verbal subgroups of some Artin groups,
Lavrent'ev (ed.), Group and metric properties of mappings. Novosibirsk: Novosibirskij
Gosudarstvennyj Universitet, 8--18 (1995), (in Russian)  (Zbl 0943.20034).
\item
    R. Gillette, J. Van Buskirk, The word problem and consequences for the
braid groups and mapping class groups of the 2--sphere, Trans.
Amer. Math. Soc., 131, N~2 (1968), 277--296.
\item
    A. I. Malcev, On isomorphic representation of infinite groups by matrix,
Matem. Sb., 8, N~3 (1940), 405--422 (in Russian).
\item
    J. Van Buskirk, Braid groups of compact 2--manifolds with elements of
finite order, Trans. Amer. Math. Soc., 122, N~1 (1966), 81--97.
\item
    M. G. Zinno, On Krammer's representation of the braid group, Math. Ann., 321, N~1 (2000),
197--211.

\end{enumerate}

\vskip 12pt

\bigskip
\bigskip
\noindent
Author address:

\bigskip

\noindent
Valerij Bardakov\\ Sobolev Institute of Mathematics,\\
Novosibirsk, 630090, Russia\\ {\tt bardakov@math.nsc.ru}

 \end{document}